\documentclass[final,hidelinks,onefignum,onetabnum]{siamart251216}
\usepackage{amsfonts,amssymb,amsmath}
\usepackage{xcolor}         
\usepackage{dsfont}
\usepackage{graphicx}
\usepackage{epstopdf}
\usepackage{algorithmic}
\usepackage{enumitem}
\ifpdf
  \DeclareGraphicsExtensions{.eps,.pdf,.png,.jpg}
\else
  \DeclareGraphicsExtensions{.eps}
\fi

\usepackage{amsmath}
\usepackage{graphicx}
\usepackage{listings}
\usepackage{subfigure}
\usepackage{subeqnarray}
\usepackage{amsfonts}
\usepackage{amssymb}
\usepackage{amsbsy}
\usepackage{amsmath}
\usepackage[mathscr]{eucal}
\usepackage{color}
\usepackage{lipsum}
\usepackage{amsfonts}
\usepackage{graphicx}
\usepackage{epstopdf}
\usepackage{algorithmic}
\usepackage{algorithm}
\usepackage{amsopn}
\usepackage{makecell}
\usepackage{multirow} 
\usepackage{tikz}
\usepackage{float}
\usetikzlibrary{arrows.meta}

\newtheorem{example}{Example}

\crefname{assumption}{Assumption}{Assumptions}
\usepackage{mathtools}

\usepackage[normalem]{ulem} 
\usepackage{dsfont} 
\usepackage{xcolor}

\newcommand{\cA}{{\mathcal{A}}}

\newcommand{\cK}{{\mathcal{K}}}

\newcommand{\cS}{{\mathcal{S}}}

\newcommand{\cX}{{\mathcal{X}}}

\newcommand{\RR}{\mathbb{R}} 

\newcommand{\bc}{\begin{center}}
\newcommand{\ec}{\end{center}}

\newcommand{\bdm}{\begin{displaymath}}
\newcommand{\edm}{\end{displaymath}}

\newcommand{\beq}{\begin{equation}}
\newcommand{\eeq}{\end{equation}}

\newcommand{\bfl}{\begin{flushleft}}
\newcommand{\efl}{\end{flushleft}}

\newcommand{\bt}{\begin{tabbing}}
\newcommand{\et}{\end{tabbing}}

\newcommand{\beqn}{\begin{eqnarray}}
\newcommand{\eeqn}{\end{eqnarray}}

\newcommand{\beqs}{\begin{align*}} 
\newcommand{\eeqs}{\end{align*}}  

\newtheorem{remark}{Remark}

\title{Constraint Qualifications and Regularity Conditions
for Nonlinear Constraint Systems}
\author{Shiji Zuo\thanks{Department of Mathematical Sciences, Rensselaer Polytechnic Institute, Troy, NY (\email{zuos@rpi.edu, xuy21@rpi.edu}).}
\and Qihang Lin\thanks{Tippie College of Business, University of Iowa, Iowa City, IA (\email{qihang-lin@uiowa.edu}).}
\and Wei Liu\thanks{IMAI, Wuhan University; Department of Applied Mathematics, Hong Kong Polytechnic University, China (\email{lwdsdqqb@gmail.com}).}\and Yangyang Xu\footnotemark[1]}
\date{\today}

\begin{document}

\maketitle

\begin{abstract}
Constraint qualifications (CQs) and regularity conditions (RCs) play a fundamental role in establishing global iteration and oracle complexity results for first-order methods on solving functional constrained optimization problems. Various conditions have been proposed and employed in the literature, with their validity established for different classes of applications. However, the use of different conditions across existing studies makes it difficult to directly compare their theoretical results and underlying assumptions. In this report, we systematically review several widely used CQs and RCs and investigate their relationships. In particular, we establish implications among these conditions and identify cases in which the implications are strict, thereby clarifying their relative strengths. Our results provide a unified perspective on existing regularity assumptions and offer guidance for selecting appropriate conditions in the complexity analysis of first-order methods.    
\end{abstract}

\section{Introduction}
In this report, we compare several constraint qualifications (CQs) and regularity conditions (RCs) used in the convergence and complexity analyses of algorithms for smooth constrained
optimization. We examine the implications among these conditions, the regions on which they are imposed, and their
behavior under slack-variable reformulations. Specifically, we consider the general constraint system
\begin{align}\label{eq:general_model}
    x\in\mathcal X,\qquad g(x)=0,\qquad h(x)\le 0,
    \tag{C0}
\end{align}
where $\mathcal X\subseteq\RR^d$ is a closed convex set, $g:\RR^d\to\RR^n$ and $h:\RR^d\to\RR^m$ are continuously differentiable, and the inequality is componentwise. Introducing a slack variable \mbox{$s\in \mathcal S:=\RR_+^m$}, we can represent \cref{eq:general_model} equivalently as 
\begin{align}
    G(x,s):=
    \begin{bmatrix}
        g(x)\\
        h(x)+s
    \end{bmatrix} = 0,
    \qquad
    (x,s)\in\Omega:=\mathcal X\times\mathcal S.
    \label{eq:slack-G-app1}
\end{align}

Every feasible point $x$ of \cref{eq:general_model} corresponds to the unique feasible pair $(x,-h(x))$ of \eqref{eq:slack-G-app1}, and every feasible pair
of \eqref{eq:slack-G-app1} yields a feasible point
of \cref{eq:general_model}.
This reformulation allows us to study mixed equality and inequality constraints within an equality-constrained framework over a closed convex set.

Two special cases will also be analyzed in this report. The first one is the equality-constrained system
\begin{align}\label{eq:model}
    x\in\mathcal X,\qquad g(x)=0,
    \tag{C1}
\end{align}
whose feasible set is given by
\begin{align}\label{eq:cK}
    \mathcal K:=\{x\in\mathcal X:g(x)=0\}.
\end{align}
The second one is the inequality-constrained system
\begin{align}\label{eq:model2}
    x\in\mathcal X,\qquad h(x)\le 0.
    \tag{C2}
\end{align}
When no equality constraints are present, the slack-variable reformulation \cref{eq:slack-G-app1} for the system \cref{eq:model2} reduces to 
\begin{align}
    G(x,s):= 
        h(x)+s  = 0,
    \qquad
    (x,s)\in\Omega:=\mathcal X\times\mathcal S,
    \label{eq:slack-G-app2}
\end{align}
where we slightly abuse notation by continuing the usage of $G(x,s)$.

CQs and related RCs play an important role in constrained optimization, not only for ensuring that local minimizers satisfy the Karush-Kuhn-Tucker (KKT) conditions~\cite{bertsekas1999nonlinear}, but also for establishing the global convergence and iteration complexity of optimization algorithms. To rule out infeasible stationary points and bound Lagrange multipliers, various conditions have been introduced in the literature, differing both in their mathematical forms and in the regions on which they are imposed.

For equality constraints over a closed convex set, Xiao et al.~\cite{xiao2026exact} impose an extended relaxed constant rank constraint qualification (rCRCQ) at every feasible point. 
Robinson's constraint qualification (RCQ) has been considered in \cite{robinson1976stability} and reduces to the classical linear independence constraint qualification (LICQ) when $\mathcal{X} = \mathbb{R}^d$.
Shi and Hong~\cite{shi2020penalty} employ RCQ in the convergence analysis of a penalty dual decomposition method; requiring RCQ to hold at every feasible point, referred to as uniform RCQ, ensures that any limit point of the iterates satisfies the regularity requirement. For equality constraints with $\cX = \RR^d_+$, Curtis et al. \cite{curtis2025proximal} assume LICQ at limit points for asymptotic convergence, and strengthen it to a sequential uniform singular-value lower bound on the active constraint Jacobian along the iterates to establish worst-case iteration complexity.
For inequality-constrained problems, Di Pillo et al.~\cite{di_pillo_truncated_2010,pillo2002AL} assume that LICQ holds at every feasible point when $\mathcal{X} = \mathbb{R}^d$. Since LICQ implies the Mangasarian-Fromovitz constraint qualification (MFCQ) in $\mathbb{R}^d$ \cite[Section 12.6]{nocedal2006numerical}, a weaker requirement is uniform MFCQ, which Boob et al.~\cite{boob2025level} impose at every feasible point over a general closed convex set $\mathcal{X}$. 
Non-singularity conditions instead bound the constraint violation by a first-order stationarity residual. 
Sahin et al.~\cite[Theorem~4.1]{sahin2019inexact} use such a condition along the iterates of an inexact augmented Lagrangian method, whereas Lin et al.~\cite[Assumption 4]{lin2022complexity} impose a global non-singularity condition over the entire set $\mathcal X$.

Different applications may satisfy one CQ or RC but fail others.  
Hence, to provide guidance for practitioners to choose or verify a certain CQ, we systematically compare these conditions in this report. Our aim is to clarify which conditions imply others within the stated setting and to identify the limits of these implications through counterexamples.

\subsection*{Notation} 

For a positive integer $m$,    let 
$\mathbb R_+^m:=\{u\in\mathbb R^m:u_i\ge0,\ i\in[m]\}$ with $[m]:=\{1,\ldots,m\}$. We use $\|\cdot\|$ for the Euclidean norm of a vector and the induced operator norm of a matrix.
The Frobenius norm is denoted by $\|\cdot\|_F$.
For a vector $u$, the vector $[u]_+$ is defined componentwisely by $([u]_+)_i:=\max\{u_i,0\}$. We use
$\sigma_r(A)$ to denote the $r$th largest singular value
of a matrix~$A$. For a continuously differentiable mapping $F:\mathbb R^p\to\mathbb R^q$, its Jacobian
$J_F(z)\in\mathbb R^{q\times p}$ has $\nabla F_i(z)^\top$ as its $i$th row.

For a nonempty set $\mathcal C$, $\operatorname{aff}(\mathcal C)$
denotes its affine hull, and
$\operatorname{dist}(z,\mathcal C):=\inf_{y\in\mathcal C}\|z-y\|$
denotes the distance from $z$ to $\mathcal C$. 
For two subsets $\mathcal C$ and $\mathcal D$ of $\mathbb R^d$, a matrix $A\in \RR^{n\times d}$ and a vector $v\in \RR^d$, we define $$A\mathcal C:=\{Au:u\in\mathcal C\},\ 
\mathcal C+\mathcal D
:=\{u+w:u\in\mathcal C, w\in\mathcal D\}, \text{ and }\mathcal C + v :=\{v+u:u\in\mathcal C\}.$$
For a nonempty closed convex set $\mathcal C$,
$\mathcal P_{\mathcal C}(z):=\arg\min_{y\in\mathcal C}\|z-y\|$
denotes the Euclidean projection of $z$ onto $\mathcal C$;
when $\mathcal C$ is a linear subspace, $\mathcal P_{\mathcal C}$
is the corresponding orthogonal projection operator.
At a point $z\in\mathcal C$,
we denote its normal and tangent  cones at $z$, respectively by $$\mathcal N_{\mathcal C}(z) := \left\{v:\langle v,y-z\rangle\le 0,\ \forall y\in\mathcal C\right\} \text{ and }\mathcal T_{\mathcal C}(z) := \left\{d:\langle d,v\rangle\le 0,\ \forall v\in\mathcal N_{\mathcal C}(z)\right\}.$$

\section{  Constraint Qualifications and Regularity Conditions}\label{sec:cq-reg}
We collect\footnote{This list will be extended with the advancement of research on solving functional constrained optimization problems.} several CQs and RCs, which are stated below. Here, ``uniform'' means that the corresponding
CQ holds at every feasible point.
\begin{enumerate}[leftmargin=*, labelindent=0pt, label=(\roman*)]
    \item \textbf{Extended  rCRCQ for \cref{eq:model}~\cite[Assumption~1.1(2)]{xiao2026exact}.}
    
There exists a positive integer $r$ such that,  for any $x\in\mathcal K$, there is $\tau_x>0$ such that\footnote{In \cite[Assumption 1.1]{xiao2026exact}, the extended rCRCQ is formulated using an open set $\mathcal{Y} \supseteq \mathcal{X}$ on which $J_g$ is locally Lipschitz continuous. Since $\mathcal{Y}$ is open, the requirement $y \in \mathcal{Y}$ in the local constant-rank condition can be omitted by shrinking $\tau_x$ if necessary. Moreover, our analysis only requires $g$ to be continuously differentiable, without assuming local Lipschitz continuity of $J_g$.}
    \begin{align}\label{eq:rcrcq}
        \begin{aligned}
        &\dim\left(\left\{J_g(y)d:d\in\mathcal E \right\}\right)
        =r,
        \qquad
        \forall y \in \{y\in \text{aff}(\cX): \|y-x\|\le\tau_x\},\\
        &\dim\left(\left\{
            J_g(y)d:
            d\in \text{lin}(\mathcal T_{\mathcal X}(y))
        \right\}\right)
        =r,\ \forall y \in \{y\in\cX: \|y-x\|\le\tau_x\},
        \end{aligned}
    \end{align} 
    where $\mathcal E :=\operatorname{aff}(\mathcal X)-x$ is the direction space of $\operatorname{aff}(\mathcal X)$
and is independent of the choice of $x\in\mathcal X$, and $\text{lin}(\mathcal T_{\mathcal X}(y)):= \mathcal T_{\mathcal X}(y)
            \cap\bigl(-\mathcal T_{\mathcal X}(y)\bigr)$.

\item \textbf{Uniform RCQ for \cref{eq:model}~\cite[Section~III]{shi2020penalty}.}

    For every feasible point $\bar x$ of \cref{eq:model}, i.e., $\bar x\in\cK$, it satisfies
\begin{align}
    J_g(\bar{x})\mathcal{T}_{\mathcal{X}}(\bar{x})=\mathbb{R}^n.
    \label{eq:rcq-equality}
\end{align}

    \item \textbf{Uniform MFCQ for \cref{eq:model2}~\cite[Assumption~3]{boob2025level}.}
    
    For every feasible point $\bar x$ of \cref{eq:model2}, there exists
    $d\in\mathcal T_{\mathcal X}(\bar x)$ such that
    \begin{align}\label{eq:mfcq}
        \langle \nabla h_i(\bar x),d\rangle<0,
        \qquad
        \forall\,i\in\mathcal A(\bar x),
    \end{align}
    where
    \(
        \mathcal A(\bar x):=
        \big\{i\in [m] :h_i(\bar x)=0\big\}
    \)
    is the active set.  

    \item \textbf{Global non-singularity condition (GNSC) for \cref{eq:general_model}~\cite[Assumption~4]{lin2022complexity}.}
    There exists a constant $\nu>0$ such that $\forall\,x\in\mathcal X$, it holds
    \begin{align}\label{eq:lin_nonsingularity}
        \nu
        \sqrt{
            \|g(x)\|^2+\|[h(x)]_+\|^2
        }
        \le
        \operatorname{dist}\!\left(
            J_g(x)^\top g(x)
            +J_h(x)^\top [h(x)]_+,
            -\mathcal N_{\mathcal X}(x)
        \right).
    \end{align} 
    
    \item \textbf{Local non-singularity condition (LNSC) for \cref{eq:model} \cite[Assumption 2]{zuo2026first}
    }
    
    There exist constants $\sigma>0$ and $B>0$ such that
    \begin{align}\label{eq:lnsc}
        \sigma\|g(x)\|
        \le
        \operatorname{dist}\!\left(
            J_g(x)^\top g(x),
            -\mathcal N_{\mathcal X}(x) 
        \right) , \quad \forall x\in \big\{x\in\mathcal X: \  \|g(x)\|\le B\big\}.
    \end{align} 
\end{enumerate}
The first three conditions are imposed at feasible points or their neighborhoods. GNSC is imposed on the entire set $\mathcal X$, whereas LNSC is imposed on a residual sublevel set.
Thus, ``local'' in LNSC refers to sufficiently small constraint violation.
For any smooth equality system $F(z)=0$ with $z\in\mathcal C$, we use the same definition of LNSC after replacing
$(g,\mathcal X)$ by $(F,\mathcal C)$. In particular, this convention applies to the slack-variable systems considered in the sequel.

\begin{figure}[H]
    \centering
    \begin{tikzpicture}[
        cq/.style={
            draw, rounded corners=1.5pt,
            line width=0.55pt,
            inner xsep=9pt, inner ysep=6pt,
            font=\small, align=center
        },
        implication/.style={-{Stealth[length=2.5mm,width=1.8mm]},
            line width=0.65pt},
        equivalence/.style={<->, line width=0.65pt},
        edge note/.style={font=\scriptsize, align=center,
            inner sep=2pt, fill=white}
    ]
    
        \node[cq] (rcq)   at (-2.7, 3.0) {Uniform RCQ};
        \node[cq] (mfcq)  at ( 2.7, 3.0) {Uniform MFCQ};
        \node[cq] (gnsc)  at (-5.1, 0)   {GNSC};
        \node[cq] (lnsc)  at ( 0,   0)   {LNSC};
        \node[cq] (rcrcq) at ( 5.1, 0)   {Extended rCRCQ};

        \draw[equivalence] (rcq.east) --
            node[edge note, above=3pt]
                {Proposition~\ref{prop:slack-rcq-mfcq}}
            (mfcq.west);

        \draw[implication] (rcq.south) --
            node[edge note, left=4pt, pos=0.53]
                {Proposition~\ref{prop:rcq-implies-local-nsc}\\
                 compact $\mathcal X_0$}
            (lnsc.north west);

        \draw[implication] (mfcq.south) --
            node[edge note, right=4pt, pos=0.53]
                {Proposition~\ref{cor:uniform-mfcq-local-nsc}\\
                 compact $\mathcal X_0$}
            (lnsc.north east);

        \draw[implication] (gnsc.east) --
            node[edge note, above=3pt]
                {Proposition~\ref{prop:assumption4}}
            node[edge note, below=3pt]
                {bounded $J_h$}
            (lnsc.west);

        \draw[implication] (rcrcq.west) --
            node[edge note, above=3pt]
                {Proposition~\ref{strass}}
            node[edge note, below=3pt]
                {compact $\mathcal X_0$}
            (lnsc.east);

    \end{tikzpicture}
    \caption{Relationships among different CQs and RCs. The implications hold under the
    additional assumptions indicated along their arrows,
    where applicable. Compact $\cX_0$ is short for the condition: there  exists $B_0 > 0$ such that $\cX_0 := \{x \in \cX : \|g(x)\| \le B_0\}$ is compact for the constraint system \cref{eq:model} and there exists $B_0 > 0$ such that $\cX_0 := \{x \in \cX : \|[h(x)]_+\| \le B_0\}$ is compact for the constraint system \cref{eq:model2}}
    \label{fig:cq-implication-diagram}
\end{figure}
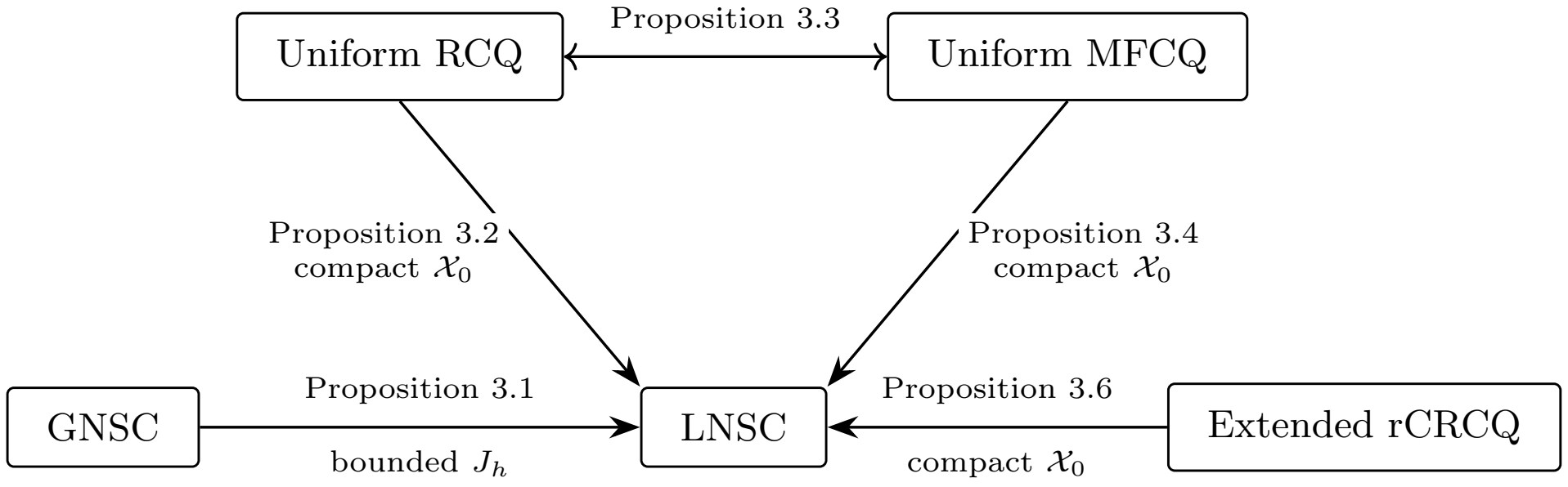

The \textbf{main results} of this report are summarized as follows.  Section~\ref{sec:gnsc-lnsc} 
shows that GNSC yields a global residual bound for the slack-variable
reformulation when $J_h$ is uniformly bounded on $\mathcal X$. Consequently, the reformulation satisfies LNSC. A counterexample shows that LNSC of the reformulation does not imply GNSC of the original system.  Section~\ref{sec:uniform-mfcq-rcq-lnsc}
establishes the equivalence between uniform MFCQ for the inequality-constrained system
\cref{eq:model2} and uniform RCQ for its slack-variable reformulation.  
 When $\mathcal X$ is compact, either condition implies LNSC, while the reverse implication fails for a counterexample. 
 In Section \ref{sec:lnsc-rcrcq}, we show that, when $\mathcal X$ is compact, extended rCRCQ implies LNSC for \cref{eq:model}; a counterexample shows that the converse does not hold in general. These implications are illustrated in \cref{fig:cq-implication-diagram}. Finally, in Section~\ref{sec:bm-clustering}, we verify these conditions for the nonnegative Burer--Monteiro relaxation arising in clustering, providing a concrete example of a structured nonconvex constraint system for which the conditions can be established directly. 

\section{Relationships among different constraint qualifications and regularity conditions}

In this section, we investigate the relationships among the aforementioned CQs and RCs listed in Section~\ref{sec:cq-reg}.
We establish several implications under their respective assumptions, and provide counterexamples
demonstrating that some of these implications are strict.

\subsection{Relationships between  GNSC and  LNSC}\label{sec:gnsc-lnsc}

We first compare the global residual condition \eqref{eq:lin_nonsingularity}, used in
\cite[Assumption~4]{lin2022complexity}, with LNSC for the slack-variable reformulation.
The following proposition shows that GNSC implies LNSC when $J_h$ is bounded on $\mathcal X$.

\begin{proposition}\label{prop:assumption4}
Consider the system \cref{eq:general_model}. Suppose that there exists a constant $L>0$ such that 
\(
    \|J_h(x)\|\le L,\, \forall\, x\in\mathcal X.
\)
If GNSC holds for \cref{eq:general_model}, then the 
reformulation \eqref{eq:slack-G-app1} satisfies LNSC. More
precisely, for every $(x,s)\in\Omega$, it holds
\begin{align}
    \sigma\|G(x,s)\|
    \le
    \operatorname{dist}\!\left(
        J_G(x,s)^\top G(x,s),
        -\mathcal N_\Omega(x,s)
    \right),
    \label{eq:ippp-implies-a1b-app}
\end{align}
where  $\sigma=
    \left[
        \frac{1}{\nu^2}
        +
        \left(1+\frac{L}{\nu}\right)^2
    \right]^{-1/2}$ with $\nu>0$ associated with GNSC in \eqref{eq:lin_nonsingularity}. 
Consequently, LNSC holds for \eqref{eq:slack-G-app1} with this constant $\sigma$ and any $B>0$.
\end{proposition}

\begin{proof}
Fixing any $(x,s)\in\Omega$, we have
\[
    \mathcal N_\Omega(x,s)
    =
    \mathcal N_{\mathcal X}(x)\times\mathcal N_{\mathcal S}(s)
    \text{ and }
    J_G(x,s)^\top G(x,s)
    =
    \begin{bmatrix}
        J_g(x)^\top g(x)+J_h(x)^\top (h(x)+s)\\
        h(x)+s
    \end{bmatrix}.
\]

By the separable structure of $\mathcal N_\Omega(x,s)$, it follows directly that
\begin{align}
    \operatorname{dist}^2\!\left(
        J_G(x,s)^\top G(x,s),
        -\mathcal N_\Omega(x,s)
    \right)
    &=
    \operatorname{dist}^2\!\left(
        J_g(x)^\top g(x)+J_h(x)^\top (h(x)+s),
        -\mathcal N_{\mathcal X}(x)
    \right) \notag\\
    &\quad +
    \operatorname{dist}^2\!\left(
        h(x)+s,-\mathcal N_{\mathcal S}(s)
    \right).
\label{eq:product-distance-app}
\end{align}

We first show that
\begin{align}
    \|h(x)+s-[h(x)]_+\|\le \operatorname{dist}\!\left(h(x)+s,-\mathcal N_{\mathcal S}(s)\right).
    \label{eq:r-p-ds-app}
\end{align}
Indeed, fix any $i\in[m]$. If $s_i=0$, then $-\mathcal N_{\mathbb R_+}(s_i)=\mathbb R_+$, and thus
\[
    |h_i(x)+s_i-[h_i(x)]_+|
    =
    |h_i(x)-[h_i(x)]_+|
    =
    \operatorname{dist}(h_i(x),\mathbb R_+)
   =\operatorname{dist}(h_i(x)+s_i,-\mathcal N_{\mathbb R_+}(s_i)).
\]
If $s_i>0$, then $-\mathcal N_{\mathbb R_+}(s_i)=\{0\}$. In this case, when $h_i(x)\le0$, we have $[h_i(x)]_+=0$ and hence $|h_i(x)+s_i-[h_i(x)]_+|=|h_i(x)+s_i|$; when $h_i(x)>0$, we have \mbox{$[h_i(x)]_+=h_i(x)$} and $0<s_i<h_i(x)+s_i$, which yields $|h_i(x)+s_i-[h_i(x)]_+|=s_i<|h_i(x)+s_i|$. Combining both cases gives
\[
    |h_i(x)+s_i-[h_i(x)]_+|
    \le
    \operatorname{dist}(h_i(x)+s_i,-\mathcal N_{\mathbb R_+}(s_i)),
    \quad \forall\, i\in[m].
\]
Summing the squares of these componentwise inequalities over $i\in[m]$ yields
\[
    \|h(x)+s-[h(x)]_+\|^2
    =
    \sum_{i=1}^m |h_i(x)+s_i-[h_i(x)]_+|^2
    \le
    \operatorname{dist}^2\!\left(h(x)+s,-\mathcal N_{\mathcal S}(s)\right),
\]
which establishes \eqref{eq:r-p-ds-app}.

Next, since the distance function to a closed set is $1$-Lipschitz and $\|J_h(x)\|\le L$, combining the GNSC condition \eqref{eq:lin_nonsingularity} with \eqref{eq:r-p-ds-app} gives
\begin{align}
    &\nu\sqrt{\|g(x)\|^2+\|[h(x)]_+\|^2}
    \le
    \operatorname{dist}\!\left(
        J_g(x)^\top g(x)+J_h(x)^\top [h(x)]_+,
        -\mathcal N_{\mathcal X}(x)
    \right) \notag\\
     \le &
    \operatorname{dist}\!\left(
        J_g(x)^\top g(x)+J_h(x)^\top (h(x)+s),
        -\mathcal N_{\mathcal X}(x)
    \right)
    +
    \left\|J_h(x)^\top\big([h(x)]_+-(h(x)+s)\big)\right\| \notag\\
     \le &
    \operatorname{dist}\!\left(
        J_g(x)^\top g(x)+J_h(x)^\top (h(x)+s),
        -\mathcal N_{\mathcal X}(x)
    \right)
    +
    L\operatorname{dist}\!\left(h(x)+s,-\mathcal N_{\mathcal S}(s)\right).
    \label{eq:ass4-bound-app}
\end{align}

Applying the triangle inequality to the block vector $G(x,s)=\big(g(x),\,h(x)+s\big)$, together with \eqref{eq:r-p-ds-app} and \eqref{eq:ass4-bound-app}, we obtain
\begin{align*}
    \|G(x,s)\| 
    &=
    \sqrt{\|g(x)\|^2+\|h(x)+s\|^2}
     \le
    \sqrt{\|g(x)\|^2+\|[h(x)]_+\|^2}+\|h(x)+s-[h(x)]_+\| \\
    &\le
    \frac{1}{\nu}\operatorname{dist}\!\left(
        J_g(x)^\top g(x)+J_h(x)^\top (h(x)+s),
        -\mathcal N_{\mathcal X}(x)
    \right) \\
    &\quad \quad\quad +
    \left(1+\frac{L}{\nu}\right)\operatorname{dist}\!\left(h(x)+s,-\mathcal N_{\mathcal S}(s)\right).
\end{align*}
Finally, applying the Cauchy--Schwarz inequality to the above inequality, and recalling~\eqref{eq:product-distance-app}, we arrive at
\begin{align*}
    \|G(x,s)\| 
    &\le
    \sqrt{
        \frac{1}{\nu^2}
        +
        \left(1+\frac{L}{\nu}\right)^2
    }
    \operatorname{dist}\!\left(
        J_G(x,s)^\top G(x,s),
        -\mathcal N_\Omega(x,s)
    \right)\\
    &=
    \frac{1}{\sigma}
    \operatorname{dist}\!\left(
        J_G(x,s)^\top G(x,s),
        -\mathcal N_\Omega(x,s)
    \right),
\end{align*}
which proves \eqref{eq:ippp-implies-a1b-app}. Since this inequality holds for every $(x,s)\in\Omega$, we conclude that LNSC holds for \eqref{eq:slack-G-app1} with the constant $\sigma$ and any $B>0$.
\end{proof}

    Conversely, the global bound \eqref{eq:ippp-implies-a1b-app} implies GNSC for the   system~\cref{eq:general_model}. Indeed, for any $x\in\mathcal X$, setting $s=[-h(x)]_+\in\mathcal S$ yields $h(x)+s=[h(x)]_+$. Since $[h(x)]_+\ge0$, $[-h(x)]_+\ge0$, and $\langle [h(x)]_+,[-h(x)]_+\rangle=0$, we have
    \[
        [h(x)]_+\in -\mathcal N_{\mathcal S}\big([-h(x)]_+\big),
    \]
    which implies $\operatorname{dist}\!\left(h(x)+[-h(x)]_+,-\mathcal N_{\mathcal S}\big([-h(x)]_+\big)\right)=0$. Consequently, evaluating \eqref{eq:product-distance-app} at $s=[-h(x)]_+$ gives
    \[
    \begin{aligned}
        &\operatorname{dist}\!\left(
            J_G\big(x,[-h(x)]_+\big)^\top G\big(x,[-h(x)]_+\big),
            -\mathcal N_\Omega\big(x,[-h(x)]_+\big)
        \right)\\
        &\qquad=
        \operatorname{dist}\!\left(
            J_g(x)^\top g(x)+J_h(x)^\top [h(x)]_+,
            -\mathcal N_{\mathcal X}(x)
        \right).
    \end{aligned}
    \]
    Moreover, since
    \[
        \left\|G\big(x,[-h(x)]_+\big)\right\|^2
        =
        \|g(x)\|^2+\|[h(x)]_+\|^2,
    \]
    evaluating \eqref{eq:ippp-implies-a1b-app} at $\big(x,[-h(x)]_+\big)$ immediately yields GNSC in \eqref{eq:lin_nonsingularity} with $\nu=\sigma$. Thus, when $J_h$ is uniformly bounded on $\mathcal X$, GNSC for \cref{eq:general_model} is equivalent to the global bound \eqref{eq:ippp-implies-a1b-app} for the slack-variable reformulation \eqref{eq:slack-G-app1}. The example below illustrates that this equivalence no longer holds under the weaker local requirement of LNSC.

\begin{example} 
In this example, we show that LNSC for the slack-variable reformulation does
not imply GNSC for the original system, even
when $\mathcal X$ is compact. Consider the system 
\begin{align}
    x\in\mathcal X:=[-1,1],
    \qquad h(x):=x(2-x)\le0.
    \label{eq:strict-example-ippp-app}
\end{align}
The function $h$ is $4$-Lipschitz continuous
on $\mathcal X$, since $$h'(x)=2-2x \text{ and }\max_{x\in\mathcal X}|h'(x)|=4.$$
Nevertheless, GNSC fails. Indeed,
$h(1)=1$, $h'(1)=0$, and
$-\mathcal N_{\mathcal X}(1)=(-\infty,0]$.
Hence
\[
    \operatorname{dist}\!\left(
        h'(1)[h(1)]_+,-\mathcal N_{\mathcal X}(1)
    \right)=0,
    \qquad [h(1)]_+=1,
\]
which rules out \eqref{eq:lin_nonsingularity}
for every $\nu>0$.

On the other hand, we show that LNSC holds for the reformulated system \cref{eq:slack-G-app2} with
\(
    B=\frac34\) and \( 
    \sigma=1.
\)
Suppose $(x,s)\in \Omega$ and  
\(    |G(x,s)|\le\frac34.
\)
If $s>0$, then $-\mathcal N_{\mathbb R_+}(s)=\{0\}$, and thus 
\[
    \operatorname{dist}\!\left(
        J_G(x,s)^\top G(x,s),-\mathcal N_\Omega(x,s)
    \right)\ge\operatorname{dist}\left(G(x,s),-\mathcal N_{\mathbb R_+}(s)\right) 
    = |G(x,s)|.
\]
If $s=0$ and $h(x)\le0$, then $G(x,s)=h(x)\le0$ and
$-\mathcal N_{\mathbb R_+}(0)=\mathbb R_+$, which again gives
\[
    \operatorname{dist}\!\left(
        J_G(x,0)^\top G(x,s),-\mathcal N_\Omega(x,0)
    \right)
    \ge
    \operatorname{dist}(G(x,s),\mathbb R_+)
    =|G(x,s)|.
\]
Finally, suppose $s=0$ and $h(x)>0$. Then $x\in(0,1]$ and
$
    G(x,s)=h(x)= 1-(1-x)^2.
$
Since $G(x,s)\le3/4$, we have $x\le1/2$, and hence
\(
    h'(x)=2(1-x)\ge1.
\) 
Moreover, it holds $x\in(0,1/2]\subset\operatorname{int}(\mathcal X)$, so
$-\mathcal N_{\mathcal X}(x)=\{0\}$. Because $G(x,s)\ge 0$ belongs to
$-\mathcal N_{\mathbb R_+}(0)=\mathbb R_+$, it follows that
\[
    \operatorname{dist}\!\left(
        J_G(x,0)^\top G(x,s),-\mathcal N_\Omega(x,0)
    \right)
    =
    |h'(x)G(x,s)|
    \ge |G(x,s)|. 
\]
Combining the above three cases yields
\[
    \|G(x,s)\|
    \le
    \operatorname{dist}\!\left(
        J_G(x,s)^\top G(x,s),
        -\mathcal N_\Omega(x,s)
    \right)
\]
whenever $(x,s)\in\Omega$ and $\|G(x,s)\|\le3/4$. 
Therefore,
LNSC holds, while GNSC fails for the instance in \eqref{eq:strict-example-ippp-app}.
\end{example}

\subsection{Relationships among uniform MFCQ, uniform RCQ and LNSC}
\label{sec:uniform-mfcq-rcq-lnsc}

We next show that uniform RCQ  
implies LNSC whenever a sufficiently small residual sublevel set is compact.
This assumption is satisfied, in particular, when $\mathcal X$ is compact.

\begin{proposition}  
\label{prop:rcq-implies-local-nsc}
Consider the constraint system in \eqref{eq:model}. Assume that there exists $B_0>0$ such that $\cX_0:=\{x\in\mathcal{X}: \|g(x)\|\le B_0\}$ is compact. In addition, assume that uniform RCQ holds, i.e.,  
\eqref{eq:rcq-equality} is satisfied at every feasible point
$\bar{x}\in \mathcal K$. Then, there exist constants $\sigma>0$ and
$B\in(0,B_0]$ such that LNSC holds.
\end{proposition}

\begin{proof}
Suppose, to the contrary,  that no such pair $(\sigma,B)$ exists.
Then for each integer $k\ge1$, we can choose a point
$x_k\in\mathcal X$ violating \eqref{eq:lnsc}
with $\sigma=1/k$ and
$B=\min\{B_0,1/k\}$.
Since \eqref{eq:lnsc} holds automatically
when $g(x)=0$, we have $g(x_k)\ne0$.
Consequently,
\[
    0<\|g(x_k)\|\le \min\left\{B_0, \frac1k\right\} \ \text{  and }\      \frac{
        \operatorname{dist}\!\left(
            J_g(x_k)^\top g(x_k),-\mathcal{N}_{\mathcal{X}}(x_k)
        \right)}
        {\|g(x_k)\|}
    \le\frac1k.
\]
Define
\(
    u_k:=\frac{g(x_k)}{\|g(x_k)\|}.
\)
Since $\mathcal{N}_{\mathcal{X}}(x_k)$ is a cone, positive homogeneity gives
\begin{align}
    \operatorname{dist}\!\left(
        J_g(x_k)^\top u_k,-\mathcal{N}_{\mathcal{X}}(x_k)
    \right)\longrightarrow0.
    \label{eq:normalized-local-nsc-failure}
\end{align}
 Because $\mathcal{X}_0$ is compact and $\{u_k\}$ lies on the unit sphere, after
passing to a subsequence, we may assume $x_k\to\bar{x}$ and $u_k\to\bar{u}$ with $\|\bar{u}\|=1$. 
Since $g(x_k)\to0$ and $g$ is continuous, we have
$g(\bar{x})=0$, so $\bar{x}\in \mathcal K$.

By \eqref{eq:normalized-local-nsc-failure}, there exists
$w_k\in -\mathcal{N}_{\mathcal{X}}(x_k)$ such that
\(
    J_g(x_k)^\top u_k-w_k\to0.
\) 
The sequence $\{J_g(x_k)^\top u_k\}$  converges to $J_g(\bar{x})^\top \bar{u}$, and thus $w_k\to \bar{w}=J_g(\bar{x})^\top \bar{u}$. 
Since the normal-cone mapping of a closed convex set has a closed
graph, it follows
\begin{align}
    J_g(\bar{x})^\top\bar{u}\in -\mathcal{N}_{\mathcal{X}}(\bar{x}).
    \label{eq:abnormal-multiplier-equality}
\end{align}

By RCQ at $\bar{x}$,
there exists $d\in \mathcal{T}_{\mathcal{X}}(\bar{x})$ such that
\(
    J_g(\bar{x})d=-\bar{u}.
\) 
Taking the inner product with $\bar{u}$ yields
\begin{align}\label{eq:neg-prod}
    \left\langle
        J_g(\bar{x})^\top\bar{u},d
    \right\rangle
    =
    -\|\bar{u}\|^2
    =
    -1.
\end{align}
On the other hand,
\eqref{eq:abnormal-multiplier-equality} and
$d\in \mathcal{T}_{\mathcal{X}}(\bar{x})$ imply, by the polarity relation between
$\mathcal{N}_{\mathcal{X}}(\bar{x})$ and $\mathcal{T}_{\mathcal{X}}(\bar{x})$,
\(
    \left\langle
        J_g(\bar{x})^\top\bar{u},d
    \right\rangle
    \ge0.
\)
This contradicts \eqref{eq:neg-prod} and  completes the proof.
\end{proof}

We now apply the RCQ to the slack-variable reformulation \cref{eq:slack-G-app2} of the inequality-constrained system~\cref{eq:model2}. 
Let $\bar{x}$ be a feasible point of \cref{eq:model2}, and define the corresponding slack variable and active index set by
\begin{align}\label{eq:s-A}
    \bar{s}:=-h(\bar{x})\in\cS,
    \qquad
    \cA(\bar{x})
    :=
    \big\{i\in [m]:h_i(\bar{x})=0\big\}.
\end{align} 
Since
$
    J_G(\bar{x},\bar{s})
    =[J_h(\bar{x}), I_m]
    \text{ and }
    \mathcal{T}_{\Omega}(\bar{x},\bar{s})
    =
    \mathcal{T}_{\mathcal{X}}(\bar{x})\times \mathcal{T}_{\cS}(\bar{s}),
$
the RCQ  for the reformulated system~\cref{eq:slack-G-app2} at $(\bar{x},\bar{s})$ takes the form
\begin{align}
    J_G(\bar{x},\bar{s})
    \mathcal{T}_{\Omega}(\bar{x},\bar{s})
    =
    \mathbb{R}^m,
    \label{eq:rcq-slack}
\end{align}
or equivalently,
\begin{align}
    J_h(\bar{x})\mathcal{T}_{\mathcal{X}}(\bar{x})
    +
    \mathcal{T}_{\cS}(\bar{s})
    =
    \mathbb{R}^m.
    \label{eq:rcq-slack-expanded}
\end{align}
The following proposition shows that RCQ for the slack-variable reformulation~\cref{eq:slack-G-app2} at $(\bar{x},\bar{s})$ is equivalent to MFCQ  for the original 
system~\cref{eq:model2} at $\bar{x}$.

\begin{proposition} 
\label{prop:slack-rcq-mfcq}
Let $\bar{x}\in \mathcal{X}$ be feasible for \eqref{eq:model2} and let
$\bar{s}:=-h(\bar{x})$. Then the following statements are equivalent:
\begin{enumerate}
    \item[(1)] the 
    RCQ   in 
    \eqref{eq:rcq-slack-expanded} holds at $(\bar{x},\bar{s})$; 

    \item[(2)] there exists
    $d\in \mathcal{T}_{\mathcal{X}}(\bar{x})$ such that the MFCQ   in \cref{eq:mfcq} holds at $\bar{x}$. 
\end{enumerate}
\end{proposition}

\begin{proof}

Since $\bar{s}_i=-h_i(\bar{x})$ for each $i\in[m]$, we have $\bar{s}_i=0$ if and only if $i\in \cA(\bar{x})$, and $\bar{s}_i>0$ for $i\notin \cA(\bar{x})$. Consequently, the tangent cone to $\mathcal{S}=\mathbb{R}_+^m$ at $\bar{s}$ is given by
\[
    \mathcal T_{\mathcal S}(\bar s)
    =
    \left\{
        q\in\mathbb R^m:
        q_i\ge0,\ \forall\, i\in\mathcal A(\bar x)
    \right\},
\]
where the components $q_i$ for $i\notin\mathcal A(\bar x)$ are unrestricted. If $\mathcal A(\bar x)=\varnothing$, then $\mathcal T_{\mathcal S}(\bar s)=\mathbb R^m$, in which case both \eqref{eq:rcq-slack-expanded} and \cref{eq:mfcq} hold trivially. Thus, without loss of generality, we assume that $\mathcal A(\bar x)\ne\varnothing$. 

Suppose first that \eqref{eq:rcq-slack-expanded} holds. Define $y\in\mathbb R^m$ by
\[
    y_i=
    \begin{cases}
        -1,&i\in\mathcal A(\bar x),\\
        0,&i\notin\mathcal A(\bar x).
    \end{cases}
\]
By \eqref{eq:rcq-slack-expanded}, there exist $d\in \mathcal{T}_{\mathcal{X}}(\bar{x})$ and $q\in \mathcal{T}_{\cS}(\bar{s})$ such that $y=J_h(\bar{x})d+q$. Since $q_i\ge0$ for all $i\in\cA(\bar{x})$, it follows that
\[
    \langle\nabla h_i(\bar{x}),d\rangle
    =
    y_i-q_i
    =
    -1-q_i
    \le -1<0,
    \qquad \forall\, i\in \cA(\bar{x}),
\]
which shows that MFCQ in \cref{eq:mfcq} holds at $\bar{x}$.

Conversely, suppose that the MFCQ holds at $\bar{x}$, so that there exists $d_0\in \mathcal{T}_{\mathcal{X}}(\bar{x})$ satisfying 
$
    \langle\nabla h_i(\bar{x}),d_0\rangle<0,$ $
  \forall\,  i\in \cA(\bar{x}).
$
For an arbitrary $y\in\mathbb{R}^m$, we can choose $t>0$ sufficiently large such that 
$
    y_i
    -
    t\langle\nabla h_i(\bar{x}),d_0\rangle
    \ge0,
    \
   \forall\, i\in \cA(\bar{x}).
$
Since $\mathcal{T}_{\mathcal{X}}(\bar{x})$ is a cone, we have $td_0\in \mathcal{T}_{\mathcal{X}}(\bar{x})$.
Let 
\(
    q:=y-J_h(\bar{x})(td_0).
\) 
Then $q_i\ge0$ for every $i\in \cA(\bar{x})$, while $q_i$ is
unrestricted for $i\notin \cA(\bar{x})$. Hence
$q\in \mathcal{T}_{\cS}(\bar{s})$, and
\(
    y=J_h(\bar{x})(td_0)+q.
\)
Since $y\in\mathbb{R}^m$ is arbitrary, this proves
\eqref{eq:rcq-slack-expanded} and completes the proof. 
\end{proof}

Each feasible point $x$ of \cref{eq:model2} corresponds to the unique feasible pair
$(x,-h(x))$ of system \cref{eq:slack-G-app2}. Thus, \cref{prop:slack-rcq-mfcq} implies that MFCQ holds at every feasible point of the
original system if and only if RCQ holds at every feasible point of its slack-variable reformulation. Combining this equivalence with
\cref{prop:rcq-implies-local-nsc} gives the following result.

\begin{proposition} 
\label{cor:uniform-mfcq-local-nsc}
Assume that there exists a positive constant $B_0$ such that the set \mbox{$\cX_0:=\{x\in\mathcal{X}: \|[h(x)]_+\|\le B_0\}$} is compact. If uniform MFCQ holds for the constraint system \cref{eq:model2}, then the slack-variable reformulation \cref{eq:slack-G-app2} satisfies LNSC. 
\end{proposition}

\begin{proof}
By \cref{prop:slack-rcq-mfcq}, for every feasible point
$\bar{x}$ of \eqref{eq:model2}, the RCQ condition in  \eqref{eq:rcq-slack-expanded} holds at $(\bar{x},\bar{s})$ with $\bar{s}=-h(\bar x)$. 
It remains to verify the compactness needed in
\cref{prop:rcq-implies-local-nsc}. 
Define
$
    \Omega_0
    :=
    \{(x,s)\in\Omega:\|G(x,s)\|\le B_0\}.
$
For any $(x,s)\in\Omega_0$, the inequality $s\ge0$ gives
$
    0\le[h(x)]_+\le[h(x)+s]_+,
$
componentwisely.
Consequently,
\[
    \|[h(x)]_+\|
    \le\|[h(x)+s]_+\|
    \le\|h(x)+s\|
    \le B_0,
\]
so $x\in\mathcal X_0$.
Since $h$ is continuous and $\mathcal X_0$
is compact, there exists $B_h\ge0$ such that
$\|h(x)\|\le B_h$ for all $x\in\mathcal X_0$.
Using $s=G(x,s)-h(x)$, we obtain
\[
    \|s\|
    \le\|G(x,s)\|+\|h(x)\|
    \le B_0+B_h.
\]
Thus $\Omega_0$ is bounded.
It is also closed, because $\Omega$ is closed
and $G$ is continuous. Therefore $\Omega_0$
is compact.
Applying \cref{prop:rcq-implies-local-nsc}
to $(G,\Omega)$ proves the claim.
\end{proof}

The example below illustrates that the converse of \cref{cor:uniform-mfcq-local-nsc} fails in general.

\begin{example}
In this example, we show that LNSC for the slack-variable reformulation does not imply uniform MFCQ for the original system or uniform RCQ for the reformulation, even when $\mathcal X$ is compact.
 Consider the constraint system
\begin{align}
    x\in\mathcal X:=[-1,1],
    \qquad
    h(x):=x^2\left(x-\frac12\right)\le0.
    \label{eq:mfcq-counterexample}
\end{align}
At the feasible point $\bar x=0$, the constraint is active with $h(\bar x)=0$ and $h'(\bar x)=0$. Thus, MFCQ fails at $\bar x$, and by \cref{prop:slack-rcq-mfcq}, RCQ also fails at $(\bar x,0)$ for the corresponding slack-variable reformulation~\cref{eq:slack-G-app2}.

On the other hand, we show that LNSC holds for \cref{eq:slack-G-app2} with $\sigma=\frac14$ and any $B>0$. Suppose $(x,s)\in\Omega$, so that $J_G(x,s)=\begin{bmatrix}h'(x) & 1\end{bmatrix}$ with $h'(x)=3x^2-x$. If $s>0$, then $-\mathcal N_{\mathbb R_+}(s)=\{0\}$, and thus
\[
    \operatorname{dist}\!\left(
        J_G(x,s)^\top G(x,s),-\mathcal N_\Omega(x,s)
    \right)
    \ge
    \operatorname{dist}\!\left(
        G(x,s),-\mathcal N_{\mathbb R_+}(s)
    \right)
    =|G(x,s)|.
\]
If $s=0$ and $h(x)\le0$, then $G(x,s)=h(x)\le0$ and $-\mathcal N_{\mathbb R_+}(0)=\mathbb R_+$, which again gives
\[
    \operatorname{dist}\!\left(
        J_G(x,0)^\top G(x,s),-\mathcal N_\Omega(x,0)
    \right)
    \ge
    \operatorname{dist}(G(x,s),\mathbb R_+)
    =|G(x,s)|.
\]
Finally, suppose $s=0$ and $h(x)>0$. Then $x\in(1/2,1]$, $G(x,s)=h(x)>0$, and hence $h'(x)=3x^2-x\ge1/4$. Moreover, for $x\in(1/2,1)\subset\operatorname{int}(\mathcal X)$, we have \mbox{$-\mathcal N_{\mathcal X}(x)=\{0\}$}, whereas at $x=1$, \mbox{$-\mathcal N_{\mathcal X}(1)=(-\infty,0]$} and \mbox{$h'(1)G(x,s)=2 h(x)>0$}; in either case, $\operatorname{dist}\!\left(h'(x)G(x,s),-\mathcal N_{\mathcal X}(x)\right)=|h'(x)G(x,s)|$. Because $G(x,s)>0$ belongs to $-\mathcal N_{\mathbb R_+}(0)=\mathbb R_+$, it follows that
\[
    \operatorname{dist}\!\left(
        J_G(x,0)^\top G(x,s),-\mathcal N_\Omega(x,0)
    \right)
    =|h'(x)G(x,s)|
    \ge\frac14|G(x,s)|.
\]
Combining the above three cases yields
\[
    \frac14\|G(x,s)\|
    \le
    \operatorname{dist}\!\left(
        J_G(x,s)^\top G(x,s),
        -\mathcal N_\Omega(x,s)
    \right)
\]
for all $(x,s)\in\Omega$. Therefore, LNSC holds with $\sigma=\frac14$ and any $B>0$, while uniform MFCQ and uniform RCQ fail for the instance in \eqref{eq:mfcq-counterexample}.
\end{example}

\subsection{Relationships between  LNSC and  extended rCRCQ}\label{sec:lnsc-rcrcq}

In this subsection, we show that extended rCRCQ implies LNSC under
the residual-sublevel compactness assumption used in
\cref{prop:rcq-implies-local-nsc}. We first establish a local estimate based on the constant-rank property.

\begin{lemma}\label{lem:local-gradient-bound}
Consider the   system in \eqref{eq:model}.
Suppose that $g$ is continuously differentiable and that
extended rCRCQ holds at $\bar x\in\mathcal K$, where $\cK$ is given in \eqref{eq:cK}. 
Let $r$ be the rank in \eqref{eq:rcrcq}.
Define
$
    \pi(\bar x)
    :=
    \sigma_r\!\left(
        \mathcal P_{\mathcal E}J_g(\bar x)^\top
    \right)>0,
$ where $\mathcal P_{\mathcal E}$ is the orthogonal projection
onto $\mathcal E$.
Then there exists a neighborhood $\mathcal U_{\bar x}$
of $\bar x$ in $\operatorname{aff}(\mathcal X)$ such that
\begin{align}
\label{eq:PE-lower}
    \left\|
        \mathcal P_{\mathcal E}J_g(y)^\top g(y)
    \right\|
    \ge
    \frac{\pi(\bar x)}{4}\|g(y)\|,
    \qquad
    \forall y\in\mathcal U_{\bar x}.
\end{align}
\end{lemma}

\begin{proof}

Let $x_0:= \arg\min_{y\in\operatorname{aff}(\mathcal X)}
\|y\|.$ 
Since $\operatorname{aff}(\mathcal X)=x_0+\mathcal E$,
we have $\mathcal P_{\mathcal E}x_0=0$ and hence
\[
y=\mathcal P_{\mathcal E}y+x_0,
\qquad
\forall y\in\operatorname{aff}(\mathcal X).
\]
Define the auxiliary mapping
$
    g^\sharp:\mathbb R^d\to\mathbb R^n$
by  $  
    g^\sharp(z)
    :=
    g(\mathcal P_{\mathcal E}z+x_0).
$
Its Jacobian matrix is 
$
    J_{g^\sharp}(z)
    =
    J_g(\mathcal P_{\mathcal E}z+x_0)
    \mathcal P_{\mathcal E}.
$
In particular, for every
$y\in\operatorname{aff}(\mathcal X)$,
\[
    g^\sharp(y)=g(y),
    \qquad
    J_{g^\sharp}(y)^\top
    =
    \mathcal P_{\mathcal E}J_g(y)^\top.
\]
Since
$\mathcal P_{\mathcal E}\bar x+x_0=\bar x$
and
$
    \|\mathcal P_{\mathcal E}z+x_0-\bar x\|
    \le \|z-\bar x\|,
$
the first rank condition in \eqref{eq:rcrcq} implies that
$J_{g^\sharp}$ has a constant rank $r$ on a neighborhood
of $\bar x$ in $\mathbb R^d$.
By the constant rank theorem~\cite[Theorem~4.12]{lee2012smooth}, we may choose a smaller
open neighborhood $\mathcal V$ of $\bar x$ such that
$
    \mathcal M:=\{g^\sharp(y):y\in\mathcal V\}
$
is an embedded $C^1$ manifold of dimension $r$
containing the origin.
Moreover,
\[
    \mathcal T_{\mathcal M}(g(y))
    =
    \operatorname{range}J_{g^\sharp}(y),
    \qquad
    \forall y\in\mathcal V\cap\operatorname{aff}(\mathcal X).
\] 

Let $Q_y$ denote the orthogonal projection onto
$\mathcal T_{\mathcal M}(g(y))$.
We next estimate the component of $g(y)$ orthogonal
to this tangent space.

Set $\mathcal L_0:=\mathcal T_{\mathcal M}(0)$.
A neighborhood of the origin in $\mathcal M$
admits a graph representation
$
    z=a+\phi(a),\ 
    \forall 
    a\in\mathcal U\subseteq\mathcal L_0,
$
where $\phi:\mathcal U\to\mathcal L_0^\perp$ is
continuously differentiable, $\phi(0)=0$,  and $J_\phi(0)=0$. For such a point $z$, the vector $a+J_\phi(a)a$ belongs to $\mathcal T_{\mathcal M}(z)$. Hence,
$
    \operatorname{dist}\bigl(
        z,\mathcal T_{\mathcal M}(z)
    \bigr)
    \le \|\phi(a)-J_\phi(a)a\|
    =o(\|a\|).
$
Since $\|z\|\ge\|a\|$, this gives
\[
    \operatorname{dist}\bigl(
        z,\mathcal T_{\mathcal M}(z)
    \bigr)
    =o(\|z\|).
\]
Applying this  to $z=g(y)$ and shrinking $\mathcal V$, we arrive at
\mbox{$
    \|g(y)-Q_y g(y)\|
    \le \frac12\|g(y)\|.
$}
Consequently,
\[
    \|Q_y g(y)\|
    \ge \frac12\|g(y)\|,
    \qquad
    \forall y\in\mathcal V\cap\operatorname{aff}(\mathcal X).
\]

On the other hand, since $J_g$ is continuous,
the continuity of the singular values gives
\[
\sigma_r\big(J_{g^\sharp}(y)^\top\big)
\ge\frac{\pi(\bar x)}2
\]
for all $y$ sufficiently close to $\bar x$.

Moreover, one has $\ker\!\left(J_{g^\sharp}(y)^\top\right)
    =
    \left(
        \operatorname{range}J_{g^\sharp}(y)
    \right)^\perp
    =
    \mathcal T_{\mathcal M}(g(y))^\perp.$ 
It follows that
\[
J_{g^\sharp}(y)^\top g(y)
=
J_{g^\sharp}(y)^\top Q_y g(y).
\]
Because $Q_y g(y)$ belongs to $\operatorname{range}J_{g^\sharp}(y)$, we obtain
\begin{align*}
&\left\|
\mathcal P_{\mathcal E}J_g(y)^\top g(y)
\right\| =\left\|
J_{g^\sharp}(y)^\top g(y)
\right\| =\left\|
J_{g^\sharp}(y)^\top Q_y g(y)
\right\|\\
&\ge
\sigma_r\big(J_{g^\sharp}(y)^\top\big)
\|Q_y g(y)\| \ge
\frac{\pi(\bar x)}2
\cdot\frac12\|g(y)\| =\frac{\pi(\bar x)}4\|g(y)\|.
\end{align*}
Taking
$\mathcal U_{\bar x}
=\mathcal V\cap\operatorname{aff}(\mathcal X)$,
after the preceding neighborhood reductions
proves the lemma.
\end{proof}

\begin{proposition}\label{strass}
Consider the   system \eqref{eq:model}.
Assume that there exists $B_0>0$ such that
$\cX_0:=\{x\in\mathcal{X}: \|g(x)\|\le B_0\}$  
is compact.
If  extended rCRCQ holds at every point of $\mathcal K$, where $\cK$ is given in \eqref{eq:cK},
then the system satisfies  LNSC in~\eqref{eq:lnsc}.
\end{proposition}

\begin{proof}
Suppose, to the contrary, that no such pair
$(\sigma,B)$ exists. Repeating the sequence construction
in the proof of \cref{prop:rcq-implies-local-nsc}, we obtain
a sequence $\{x_k\}\subset\mathcal X_0$ with
$0<\|g(x_k)\|\to0$ and unit vectors
$
    u_k:=\frac{g(x_k)}{\|g(x_k)\|}
$
such that, after passing to a subsequence,
$x_k\to\bar x\in\mathcal K$ and $u_k\to\bar u$ with
$\|\bar u\|=1$. As in
\eqref{eq:abnormal-multiplier-equality}, these limits satisfy
$
    J_g(\bar x)^\top\bar u
    \in-\mathcal N_{\mathcal X}(\bar x).
$
Consequently,
$\bar w:=-J_g(\bar x)^\top\bar u
\in\mathcal N_{\mathcal X}(\bar x)$ satisfies
\begin{align}
    J_g(\bar x)^\top\bar u+\bar w=0.
    \label{eq:Jg-u-w-zero}
\end{align}

Since $x_k\in\mathcal X\subseteq\operatorname{aff}(\mathcal X)$
and $x_k\to\bar x$, the estimate \eqref{eq:PE-lower}
in \cref{lem:local-gradient-bound} applies to $y=x_k$
for all sufficiently large $k$. Dividing this estimate
by $\|g(x_k)\|$ gives
\[
    \|\mathcal P_{\mathcal E}J_g(x_k)^\top u_k\|
    \ge \frac{\pi(\bar x)}{4}.
\]
Passing to the limit and using the continuity of $J_g$,
we obtain
\begin{align}
    \|\mathcal P_{\mathcal E}J_g(\bar x)^\top\bar u\|
    \ge
    \frac{\pi(\bar x)}{4}
    >0.
    \label{eq:PE-positive}
\end{align}

To derive a contradiction, set
\[
    \mathcal L
    :=
    \operatorname{lin}\mathcal T_{\mathcal X}(\bar x)
    =
    \mathcal T_{\mathcal X}(\bar x)
    \cap
    \bigl(-\mathcal T_{\mathcal X}(\bar x)\bigr).
\]
Since $\mathcal L\subseteq\mathcal E$, we have
$J_g(\bar x)\mathcal L
\subseteq J_g(\bar x)\mathcal E$.
The two rank conditions in \eqref{eq:rcrcq},
evaluated at $\bar x$, give
$
    \dim\bigl(J_g(\bar x)\mathcal L\bigr)
    =
    \dim\bigl(J_g(\bar x)\mathcal E\bigr)
    =r.
$
Thus,
$
    J_g(\bar x)\mathcal L
    =
    J_g(\bar x)\mathcal E.
$
Moreover, $\bar w\in\mathcal N_{\mathcal X}(\bar x)$
is orthogonal to $\mathcal L$, because both $d$ and $-d$
belong to $\mathcal T_{\mathcal X}(\bar x)$
for every $d\in\mathcal L$. It follows from
\eqref{eq:Jg-u-w-zero} that
\[
    \langle\bar u,J_g(\bar x)d\rangle
    =
    -\langle\bar w,d\rangle
    =
    0,
    \qquad d\in\mathcal L.
\]
Therefore, $\bar u$ is orthogonal to
$J_g(\bar x)\mathcal L=J_g(\bar x)\mathcal E$, and hence
$
    \mathcal P_{\mathcal E}J_g(\bar x)^\top\bar u=0.
$
This contradicts \eqref{eq:PE-positive} and completes
the proof.
\end{proof}

The example below illustrates that 
the converse of \cref{strass} fails in general.

\begin{example}
In this example,  we show that LNSC does not imply the extended rCRCQ, even when $\mathcal X$ is compact. 
Consider the constraint system
\begin{align}
    x\in\mathcal X:=[0,2]^2,
    \qquad
    g(x):=(x_1+1)^2+(x_2-1)^2-2=0.
    \label{eq:rcrcq-counterexample}
\end{align}
At the feasible point $\bar x=(0,0)$, we have
\[
    J_g(\bar x)=\begin{bmatrix}2 & -2\end{bmatrix},
    \qquad
    \mathcal E=\mathbb R^2,
    \qquad
    \mathcal T_{\mathcal X}(\bar x)=\mathbb R_+^2.
\]
Since $\operatorname{lin}\bigl(\mathcal T_{\mathcal X}(\bar x)\bigr)=\{0\}$, it follows that
\[
    \dim\bigl(J_g(\bar x)\mathcal E\bigr)=1
    \ne
    0=\dim\Bigl(J_g(\bar x)\operatorname{lin}\bigl(\mathcal T_{\mathcal X}(\bar x)\bigr)\Bigr),
\]
and thus  extended rCRCQ  fails at $\bar x$.

On the other hand, we show that LNSC holds with $\sigma=2$ and any $B>0$. Fix any $x\in\mathcal X$, and note that the first component of $J_g(x)^\top g(x)$ is $2(x_1+1)g(x)$. If $0<x_1<2$, then $-\mathcal N_{[0,2]}(x_1)=\{0\}$. If $x_1=0$, then $x_2\in[0,2]$ implies $g(x) \le0$ and hence $2(x_1+1)g(x)\le0$, while $-\mathcal N_{[0,2]}(0)=\mathbb R_+$. If $x_1=2$, then $g(x)=7+(x_2-1)^2>0$ and hence $2(x_1+1)g(x)>0$, while $-\mathcal N_{[0,2]}(2)=\mathbb R_-$. 
Consequently, in all three cases, we have
\[
    \operatorname{dist}\!\left(
        2(x_1+1)g(x),-\mathcal N_{[0,2]}(x_1)
    \right)
    =
    2(x_1+1)|g(x)|.
\]
Using the product structure $\mathcal N_{\mathcal X}(x)=\mathcal N_{[0,2]}(x_1)\times\mathcal N_{[0,2]}(x_2)$ and $x_1\ge0$, we obtain
\[
\begin{aligned}
    \operatorname{dist}\!\left(
        J_g(x)^\top g(x),
        -\mathcal N_{\mathcal X}(x)
    \right)
    &\ge
    \operatorname{dist}\!\left(
        2(x_1+1)g(x),-\mathcal N_{[0,2]}(x_1)
    \right)\\
    &=2(x_1+1)|g(x)|\ge2|g(x)|
\end{aligned}
\]
for all $x\in\mathcal X$. Therefore, LNSC holds with $\sigma=2$ and any $B>0$, while  extended rCRCQ fails for the instance in \eqref{eq:rcrcq-counterexample}. Together with \cref{strass}, this shows that extended rCRCQ is strictly stronger than LNSC when $\mathcal X$ is compact.
\end{example}

\section{A nonnegative Burer--Monteiro formulation for clustering problems}
\label{sec:bm-clustering}

Let $z_1,\ldots,z_n$ be data points in a Euclidean space,
and let $D\in\mathbb R^{n\times n}$ be their squared
Euclidean distance matrix with
$
    D_{ij}=\|z_i-z_j\|^2.
$
For a prescribed number $s$ of clusters, the semidefinite
relaxation of Peng and Wei \cite{peng2007approximating}
can be written, up to a positive scaling of the objective, as
\[
    \min_{Y\in\mathbb R^{n\times n}}
    \operatorname{tr}(DY),
    \quad\text{s.t.}\quad
    Ye=e,\quad
    \operatorname{tr}(Y)=s,\quad
    Y\succeq0,\quad Y\ge0,
\]
where $e\in\mathbb R^n$ is the vector of all ones.

Following \cite[Section~6.1]{sahin2019inexact},
we study a nonnegative factorized model based on the
Burer--Monteiro factorization $Y=VV^\top$
\cite{burer2003nonlinear}.
The model replaces the trace equality by
$\|V\|_F^2\le s$ and imposes $V\ge0$, which ensures
$VV^\top\ge0$.
This construction combines a relaxation of the trace
constraint with a restriction to nonnegative factors
having $r$ columns.
The resulting nonconvex problem is
\begin{align}
\label{eq:bm-problem}
    \min_{V\in\mathbb R^{n\times r}}
    \operatorname{tr}(DVV^\top),
    \quad\text{s.t.}\quad
    VV^\top e=e,\quad
    \|V\|_F^2\le s,\quad V\ge0.
\end{align}

In what follows, we restrict our attention to the nontrivial clustering
setting 
\begin{align}
\label{eq:bm-parameter}
    2\le s\le\min\{n,r\}.
\end{align}
The condition $r\ge s$ allows the normalized indicator
factor of any partition into $s$ nonempty clusters
to be represented in $\mathbb R^{n\times r}$
by appending zero columns;
see also \mbox{\cite[Section~2.1.2]{kulis2007fast}}.
The inequality $s>1$ will be used in the   proof below.

Equivalently,
\cref{eq:bm-problem} can be rewritten as
\begin{align}
\label{eq:bm-equality-form}
    \min_{V\in\mathcal X}
    f(V):=\operatorname{tr}(DVV^\top)
    \quad
    \text{s.t.}\quad
    g(V):=VV^\top e-e=0,
\end{align}
where
\begin{align}
\label{eq:bm-X}
    \mathcal X
    :=
    \left\{
        V\in\mathbb R^{n\times r}:
        V\ge0,\ \|V\|_F^2\le s
    \right\}.
\end{align}

For any $V\in\mathbb R^{n\times r}$, let
$J_g(V):\mathbb R^{n\times r}\to\mathbb R^n$
denote the Jacobian of $g$ at $V$, viewed as a
linear operator. For any $H\in\mathbb R^{n\times r}$, we have
\[
    J_g(V)[H]=HV^\top e+VH^\top e.
\]
For any subset $\mathcal D\subseteq\mathbb R^{n\times r}$,
we write $J_g(V)\mathcal D
    :=\{J_g(V)[H]:H\in\mathcal D\}$. 
The adjoint operator
$J_g(V)^*:\mathbb R^n\to\mathbb R^{n\times r}$
is defined by $\langle J_g(V)[H],y\rangle
    =\langle H,J_g(V)^*y\rangle$, 
for all $H\in\mathbb R^{n\times r}$ and
$y\in\mathbb R^n$. Indeed, for any $H\in\mathbb R^{n\times r}$ and
$y\in\mathbb R^n$, we have
\begin{align*}
    \langle J_g(V)[H],y\rangle
    &= y^\top(HV^\top e+VH^\top e) = \langle H,(ye^\top+ey^\top)V\rangle.
\end{align*}
Therefore, the adjoint operator is given by
\begin{align*}
    J_g(V)^*y=(ye^\top+ey^\top)V.
\end{align*}
We first establish uniform RCQ for the constraint system in
\cref{eq:bm-equality-form}.
Compactness of $\mathcal X$ then yields LNSC.

\begin{proposition}
\label{prop:bm-rcq}
Under the condition in \cref{eq:bm-parameter}, the constraint system in \cref{eq:bm-equality-form} satisfies uniform RCQ, i.e., $J_g(\bar V)\mathcal T_{\mathcal X}(\bar V)=\mathbb R^n$ for every feasible point $\bar V$,and consequently satisfies LNSC.
\end{proposition}

\begin{proof}
Let $\bar V\in\mathcal X$ be an arbitrary feasible point of problem \cref{eq:bm-equality-form}. Feasibility of $\bar V$ gives $\bar V\bar V^\top e=e$, and hence
\begin{align}
\label{eq:q-norm}
    \|\bar V^\top e\|^2
    =
    e^\top\bar V\bar V^\top e
    =
    e^\top e
    =
    n.
\end{align}
Note that $J_g(\bar V)[H]=H\bar V^\top e+\bar VH^\top e$ for any $H\in\mathbb R^{n\times r}$. For any $y\in\mathbb R^n$, since $I_n+\bar V\bar V^\top\succ0$, setting $a:=(I_n+\bar V\bar V^\top)^{-1}y$ and $H_0:=\operatorname{Diag}(a)\bar V$ yields
\begin{align}
\label{eq:bm-jacobian-surjective}
    J_g(\bar V)[H_0]
    =
    \operatorname{Diag}(a)\bar V\bar V^\top e+\bar V\bar V^\top a
    =
    (I_n+\bar V\bar V^\top)a
    =
    y,
\end{align}
which shows that $J_g(\bar V):\mathbb R^{n\times r}\to\mathbb R^n$ is surjective.
 
Next, we construct a direction $W\in\ker J_g(\bar V)$ that points strictly into the active constraints of $\mathcal X$ at $\bar V$. Let $\mathcal J_0:=\{j\in[r]:(\bar V^\top e)_j=0\}$, and note that $\bar V\ge0$ implies $\bar V_{:j}=0$ for all $j\in\mathcal J_0$. Define $Z\in\mathbb R^{n\times r}$ by $Z_{ij}:=1$ if $j\in\mathcal J_0$ and $Z_{ij}:=0$ otherwise, and set
\[
    W:=\frac{1}{n}ee^\top\bar V-\bar V+Z.
\]
Since $Z_{:j}=0$ for $j\notin\mathcal J_0$ while $(\bar V^\top e)_j=0$ and $\bar V_{:j}=0$ for $j\in\mathcal J_0$, we have $Z\bar V^\top e=0$ and $\bar VZ^\top e=0$. Combining this with \cref{eq:q-norm}, $e^\top e=n$, and $\bar V\bar V^\top e=e$ gives
\begin{align*}
    W\bar V^\top e
    &=
    \frac{\|\bar V^\top e\|^2}{n}e-\bar V\bar V^\top e+Z\bar V^\top e
    =
    0,
    \\
    \bar VW^\top e
    &=
    \frac{e^\top e}{n}\bar V\bar V^\top e-\bar V\bar V^\top e+\bar VZ^\top e
    =
    0, 
\end{align*}
and consequently
\begin{align}
\label{eq:JD-zero}
    J_g(\bar V)[W]=W\bar V^\top e+\bar VW^\top e=0.
\end{align}
Moreover, whenever $\bar V_{ij}=0$, we have $W_{ij}=(\bar V^\top e)_j/n+Z_{ij}=(\bar V^\top e)_j/n>0$ if $j\notin\mathcal J_0$ and $W_{ij}=1>0$ if $j\in\mathcal J_0$; thus,
\begin{align}
\label{eq:active-nonnegative-direction}
    W_{ij}>0 \quad\text{whenever }\bar V_{ij}=0.
\end{align}
Since $\langle\bar V,ee^\top\bar V\rangle=\|\bar V^\top e\|^2=n$ and $\langle\bar V,Z\rangle=0$, we obtain {$\langle\bar V,W\rangle=1-\|\bar V\|_F^2$}. If the ball constraint is active (i.e., $\|\bar V\|_F^2=s$), then \cref{eq:bm-parameter} gives $\langle\bar V,W\rangle=1-s<0$. 

Now define $H:=H_0+tW$ with $t>0$ chosen sufficiently large such that, by \cref{eq:active-nonnegative-direction}, $H_{ij}>0$ whenever $\bar V_{ij}=0$, and $\langle\bar V,H\rangle<0$ if $\|\bar V\|_F^2=s$. Then, for all sufficiently small $\alpha>0$, we have $$\bar V+\alpha H\ge0 \text{ and }\|\bar V+\alpha H\|_F^2=s+2\alpha\langle\bar V,H\rangle+\alpha^2\|H\|_F^2\le s$$ when $\|\bar V\|_F^2=s$ (and similarly $\|\bar V+\alpha H\|_F^2\le s$ by continuity when $\|\bar V\|_F^2<s$). Hence $\bar V+\alpha H\in\mathcal X$ for all small $\alpha>0$, which implies $H\in\mathcal T_{\mathcal X}(\bar V)$. By \cref{eq:bm-jacobian-surjective,eq:JD-zero}, we have $$J_g(\bar V)[H]=J_g(\bar V)[H_0]+tJ_g(\bar V)[W]=y.$$ Since $y\in\mathbb R^n$ and the feasible point $\bar V$ are arbitrary, uniform RCQ holds. Finally, because $\mathcal X$ is compact and $g$ is continuously differentiable, $\{V\in\mathcal X:\|g(V)\|\le B_0\}$ is compact for any $B_0>0$, and thus  LNSC follows directly from \cref{prop:rcq-implies-local-nsc}.
The proof is then completed.
\end{proof}

The following two remarks explain that, under the condition in \cref{eq:bm-parameter}, the constraint system in \cref{eq:bm-equality-form} violates both GNSC and extended rCRCQ, even though it satisfies uniform RCQ and LNSC by \cref{prop:bm-rcq}.

\begin{remark}
The constraint system in \cref{eq:bm-equality-form} does not satisfy GNSC. Indeed, at $V=0\in\mathcal X$, we have $g(0)=-e$ and $J_g(0)=0$. Since $0\in\mathcal N_{\mathcal X}(0)$, it follows that
\[
    \operatorname{dist}\!\left(
        J_g(0)^*g(0),
        -\mathcal N_{\mathcal X}(0)
    \right)
    =
    0,
    \qquad
    \|g(0)\|=\sqrt n>0.
\]
Therefore, no constant $\nu>0$ can satisfy GNSC in \cref{eq:lin_nonsingularity} on $\mathcal X$.
\end{remark}

\begin{remark}\label{exp:bmrcrcq}
Under the condition in \cref{eq:bm-parameter}, the constraint system in \cref{eq:bm-equality-form} does not satisfy the extended rCRCQ. Since $s\in\mathbb N_+$ satisfies $s\le\min\{n,r\}$ by \cref{eq:bm-parameter}, we can partition $[n]$ into $s$ nonempty disjoint subsets $C_1,\ldots,C_s$ with cardinalities $m_j:=|C_j|$ for $j\in[s]$, and define $\bar V\in\mathbb R^{n\times r}$ by
\begin{align}
\label{eq:cluster-indicator-V}
    \bar V_{ij}
    :=
    \begin{cases}
        \dfrac{1}{\sqrt{m_j}},
        & \text{if }i\in C_j\text{ for some }j\in[s],\\[1mm]
        0,
        & \text{otherwise}.
    \end{cases}
\end{align}
Then $(\bar V^\top e)_j=\sqrt{m_j}$ for $j\in[s]$ and $(\bar V^\top e)_j=0$ for $j>s$. For each $i\in[n]$, letting $j\in[s]$ be the unique index such that $i\in C_j$, we have $\bar V_{ik}=0$ for all $k\ne j$, and hence
\[
    (\bar V\bar V^\top e)_i
    =\sum_{k=1}^r \bar V_{ik}(\bar V^\top e)_k
    =\bar V_{ij}(\bar V^\top e)_j
    =\frac{1}{\sqrt{m_j}}\sqrt{m_j}=1.
\]
Thus $\bar V\bar V^\top e=e$, which means $g(\bar V)=0$. Moreover, since it holds $\bar V\ge0$ and $\|\bar V\|_F^2=\sum_{j=1}^s m_j(1/m_j)=s$, the point $\bar V$ is feasible for \cref{eq:bm-equality-form} and the ball constraint is active at $\bar V$.

Next, notice that the matrix with all entries equal to $\varepsilon\in\bigl(0,\sqrt{s/(nr)}\bigr)$ lies in $\operatorname{int}(\mathcal X)$, which implies $\operatorname{aff}(\mathcal X)=\mathbb R^{n\times r}$. Since $J_g(\bar V):\mathbb R^{n\times r}\to\mathbb R^n$ is surjective by~\cref{eq:bm-jacobian-surjective}, we have
\begin{align}
\label{eq:rcrcq-aff-dim}
    \dim\Bigl(J_g(\bar V)\bigl(\operatorname{aff}(\mathcal X)-\bar V\bigr)\Bigr)
    =
    n.
\end{align}
On the other hand, since $\mathcal X=\mathbb R_+^{n\times r}\cap\{V\in\mathbb R^{n\times r}:\|V\|_F^2\le s\}$ and $\|\bar V\|_F^2=s$, any tangent direction $H\in\mathcal T_{\mathcal X}(\bar V)$ satisfies $H_{ij}\ge0$ whenever $\bar V_{ij}=0$ as well as $\langle\bar V,H\rangle\le0$. Applying both conditions to $\pm H$ for any $H\in\operatorname{lin}\bigl(\mathcal T_{\mathcal X}(\bar V)\bigr)=\mathcal T_{\mathcal X}(\bar V)\cap\bigl(-\mathcal T_{\mathcal X}(\bar V)\bigr)$ yields
\[
    \operatorname{lin}\bigl(\mathcal T_{\mathcal X}(\bar V)\bigr)
    \subseteq
    \bigl\{H\in\mathcal W:\langle\bar V,H\rangle=0\bigr\},
    \quad\text{where}
    \]\[
    \mathcal W:=\bigl\{H\in\mathbb R^{n\times r}:H_{ij}=0\text{ whenever }\bar V_{ij}=0\bigr\}.
\]
Because $\bar V$ has exactly one positive entry in each row, we have $\dim\mathcal W=n$. Moreover, since $\bar V\in\mathcal W$ and $\langle\bar V,\bar V\rangle=s>0$, the linear functional $H\mapsto\langle\bar V,H\rangle$ is nonzero on $\mathcal W$, so its kernel in $\mathcal W$ has dimension $n-1$. Consequently,
\begin{align}
\label{eq:rcrcq-tangent-dim}
    \dim\Bigl(J_g(\bar V)\operatorname{lin}\bigl(\mathcal T_{\mathcal X}(\bar V)\bigr)\Bigr)
    \le
    \dim\operatorname{lin}\bigl(\mathcal T_{\mathcal X}(\bar V)\bigr)
    \le n-1.
\end{align}
Comparing \cref{eq:rcrcq-aff-dim,eq:rcrcq-tangent-dim} shows that the rank condition in \eqref{eq:rcrcq} is violated at $\bar V$. Therefore,  extended rCRCQ fails at the feasible point $\bar V$, even though uniform RCQ and LNSC hold for the constraint system in \cref{eq:bm-equality-form} by \cref{prop:bm-rcq}.
\end{remark}

\section{Conclusion}

 We have established the relationships among various constraint qualifications and regularity conditions. We have shown that conditions such as uniform RCQ, extended rCRCQ, uniform MFCQ, and GNSC generally imply LNSC under mild compactness or boundedness assumptions. However, as demonstrated by our counterexamples and the nonnegative Burer-Monteiro clustering model, the converse implications do not always hold. Overall, this comparison offers practical guidance for researchers to choose or verify appropriate conditions in algorithmic analysis.

\bibliographystyle{siamplain}
\bibliography{optim}

\section*{Clarification of AI Assistance}
During the preparation of this report, the authors utilized ChatGPT (GPT-6) to assist with literature searches, mathematical proof refinement, and \LaTeX{} maintenance. In particular, GPT-6 contributed to formulating \cref{lem:local-gradient-bound}, the proof of \cref{prop:bm-rcq}, and the counterexample in \cref{exp:bmrcrcq}. These arguments were subsequently verified, revised, and organized by the authors, who take full responsibility for the content of this report.
\end{document}